\documentstyle[12pt]{article}
\begin{document}
\newcommand{\R}{\bf R}
\newcommand{\F}{\bf F}
\newcommand{\wt}{\widetilde}
\vsize=217mm
\hfill
\vskip3.0 truecm
\thispagestyle{empty}
\centerline{\Large Discrete Versions of the Beckman-Quarles Theorem}
\vskip 1.3truecm
\centerline{\large Apoloniusz Tyszka}
\vskip 1.3truecm
\rightline{}
\vskip 0.3 truecm
\centerline{to appear in Aequationes Mathematicae}
\vskip 0.7 truecm
{\bf Abstract}.
Let $\F \subseteq \R$ denote the field of numbers which are
constructible by means of ruler and compass. We prove that:
{\large (1)}
if $x,y\in {\R}^{n}$ $(n>1)$ and $|x-y|$ is an algebraic number
then there exists a finite set $S_{xy}\subseteq {\R}^{n}$ containing
$x$ and $y$ such that each map from $S_{xy}$ to ${\R}^{n}$ 
preserving all unit distances preserves the distance between
$x$ and $y$; if $x,y\in {\F}^{n}$ then we can choose
$S_{xy}\subseteq {\F}^{n}$,
{\large (2)}
only algebraic distances $|x-y|$ have the property from item (1),
{\large (3)}
if $X_{1},X_{2},...,X_{m}\in {\R}^{n}$ $(n>1)$ lie on some affine
hyperplane then there exists a finite set $L(X_{1},X_{2},...,X_{m})
\subseteq {\R}^{n}$ containing $X_{1},X_{2},...,X_{m}$ such that each
map from $L(X_{1},X_{2},...,X_{m})$ to ${\R}^{n}$ preserving all unit
distances preserves the property that $X_{1},X_{2},...,X_{m}$ lie on
some affine hyperplane,
{\large (4)}
if $J,K,L,M$ $\in $ ${\R}^{n}$ $(n>1)$ and $|JK|=|LM|$ $(|JK|<|LM|)$
then there exists a finite set $C_{JKLM}\subseteq {\R}^{n}$ containing
$J,K,L,M$ such that any map $f:C_{JKLM}\rightarrow {\R}^{n}$ that
preserves unit distance satisfies $|f(J)f(K)|=|f(L)f(M)|$
($|f(J)f(K)|<|f(L)f(M)|$).
\vskip 0.3truecm
Let ${\F}\subseteq {\R}$ denote the field of numbers which are
constructible by means of ruler and compass. Theorem 1 may be
viewed as a discrete form of the classical Beckman-Quarles theorem,
which states that any map from ${\R}^{n}$ to ${\R}^{n}$
($2\leq n<\infty $) preserving unit distances is an isometry, see
[1], [2] and [5]. Theorem 1 was proved in [10] in the special case
where $n=2$ and the distance $|x-y|$ is constructible by means of
ruler and compass.
\vskip 0.1truecm
{\bf Theorem 1}. If $x,y\in {\R}^{n}$ ($n>1$) and $|x-y|$ is an
algebraic number then there exists a finite set
$S_{xy}\subseteq {\R}^{n}$ containing $x$ and $y$ such that each
map from $S_{xy}$ to ${\R}^{n}$ preserving all unit distances
preserves the distance between $x$ and $y$.
\\
{\it Proof.} The proof falls naturally into two sections.
\\
{\bf 1. Elementary facts}

Let us denote by $D_{n}$ the set of all non-negative numbers $d$
with the following property:

If $x,y\in {\R}^{n}$ and $|x-y|=d$ then there exists a finite set
$S_{xy}\subseteq {\R}^{n}$ such that $x,y\in S_{xy}$ and any map
$f:S_{xy}\rightarrow {\R}^{n}$ that preserves unit distance also
preserves the distance between $x$ and $y$.

Obviously $0,1\in D_{n}$. We first prove that if $d\in D_{n}$ then
$\sqrt{2+2/n}\cdot d\in D_{n}$. Let us assume that $d>0$,
$x,y\in {\R}^{n}$, $|x-y|=\sqrt{2+2/n}\cdot d$.
Using the notation of Figure 1 we show that
\vskip 0.2truecm
\centerline{$S_{xy}:=
\bigcup \{S_{ab}:a,b\in \{x,y,\wt{y},p_{1},p_{2},...,p_{n},
\wt{p}_{1},\wt{p}_{2},...,\wt{p}_{n}\},|a-b|=d\}$}
\vskip 0.2truecm
satisfies the condition of Theorem 1. 
Figure 1 shows the case $n=2$, but equations below Figure 1
describe the general case $n\geq 2$; $z$ denotes
the centre of the $(n-1)$-dimensional regular simplex
$p_{1}p_{2}...p_{n}$.
\\
\special{em: graph Figure1.gif}
\vskip 7.7truecm
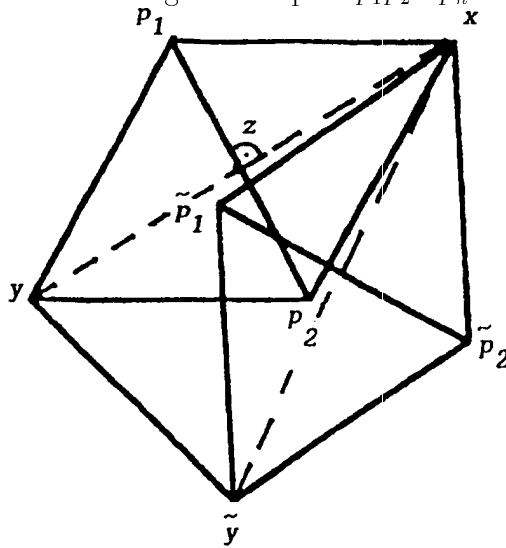
\centerline{Figure 1}
\centerline{$1\leq i<j\leq n$}
$|y-\wt{y}|=d$, $|x-p_{i}|=|y-p_{i}|=|p_{i}-p_{j}|=d=|x-\wt{p}_{i}|=
|\wt{y}-\wt{p}_{i}|=|\wt{p}_{i}-\wt{p}_{j}|$
\centerline{
$|x-\wt{y}|=|x-y|=2\cdot |x-z|=2\cdot \sqrt{\frac{n+1}{2n}}
\cdot d=\sqrt{2+2/n}\cdot d$}
\\
\vskip 1.0truecm
Let us assume that $f:S_{xy}$ $\rightarrow$ ${\R}^{n}$
preserves the distance $1$. Since
\vskip 0.3truecm
\centerline{$S_{xy}\supseteq S_{y\wt{y}}
\cup
\bigcup_{i=1}^{n}S_{xp_{i}}
\cup
\bigcup_ {i=1}^{n}S_{yp_{i}}
\cup
\bigcup_{1 \leq i<j \leq n} S_{p_{i}p_{j}}$}
\vskip 0.3 truecm
we conclude that $f$ preserves the distances between $y$ and $\wt{y}$,
$x$ and $p_{i}$ ($1\leq i \leq n$), $y$ and $p_{i}$ ($1\leq i\leq n$),
and all distances between $p_{i}$ and $p_{j}$ ($1\leq i<j\leq n$).
Hence $|f(y)-f(\wt{y})|=d$ and $|f(x)-f(y)|$ is either $0$ or
$\sqrt{2+2/n}\cdot d$. Analogously we have that $|f(x)-f(\wt{y})|$
is either $0$ or $\sqrt{2+2/n}\cdot d$. Thus $f(x)\neq f(y)$, so
$|f(x)-f(y)|=\sqrt{2+2/n}\cdot d$ which completes the proof that
$\sqrt{2+2/n}\cdot d\in D_{n}$.
\vskip 0.3truecm
Hence, if $d\in D_{n}$ then
$(2+2/n)\cdot d=\sqrt{2+2/n}\cdot (\sqrt{2+2/n}\cdot d)\in D_{n}$.
\vskip 0.3truecm
We next prove that if $x,y\in {\R}^{n}$, $d\in D_{n}$
and $|x-y|=(2/n)\cdot d$ then there exists a finite set
$Z_{xy}\subseteq {\R}^{n}$ containing $x$ and $y$
such that any map $f:Z_{xy} \rightarrow {\R}^{n}$ that preserves
unit distance satisfies $|f(x)-f(y)|\leq |x-y|$; this result is adapted
from [5]. It is obvious if $n=2$, therefore we assume that $n>2$ and
$d>0$.
Let us see at Figure 2 below, $z$ denotes the centre of the
$(n-1)$-dimensional regular simplex $p_{1}p_{2}...p_{n}$.
Figure 2 shows the case $n=3$, but equations below Figure 2 describe
the general case $n\geq 3$.
\\
\special{em: graph figure2.gif}
\vskip 5.2truecm
\centerline{Figure 2}
\centerline{$1 \leq i<j \leq n$}
\centerline{$|x-p_{i}|=|y-p_{i}|=d, \hspace{0.20cm}
|p_{i}-p_{j}|=\sqrt{2+2/n} \cdot d, \hspace{0.20cm}
|z-p_{i}|=\sqrt{1-1/n^{2}} \cdot d$}
\centerline{$|x-y|=2 \cdot |x-z|=2 \cdot
\sqrt{|x-p_{i}|^{2}-|z-p_{i}|^2}=
2 \cdot \sqrt{d^{2}-(1-1/n^{2}) \cdot d^{2}}=(2/n) \cdot d$}
\vskip 0.4truecm
Let us define:
\vskip 0.2truecm
\centerline{$Z_{xy}:=\bigcup_{1\le i<j \leq n} S_{p_{i}p_{j}}
\cup
\bigcup_{i=1}^{n} S_{xp_{i}}
\cup
\bigcup_{i=1}^{n} S_{yp_{i}}$}
\vskip 0.2truecm
If $f:Z_{xy}\rightarrow {\R}^{n}$ preserves the distance $1$ then
$|f(x)-f(y)|=|x-y|=(2/n)\cdot d$ or $|f(x)-f(y)|=0$, hence
$|f(x)-f(y)|\leq|x-y|$.
\vskip 0.2truecm
If $d\in D_{n}$, then $2\cdot d\in D_{n}$ (see Figure 3).
\\
\special{em: graph figure3.gif}
\vskip 1.2truecm
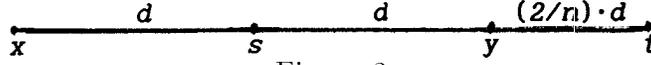
\centerline{Figure 3}
\centerline{$|x-y|=2 \cdot d$}
\centerline{$S_{xy}=S_{xs}\cup S_{sy}\cup Z_{yt}\cup S_{xt}$}
\vskip 0.5truecm
 From Figure 4 it is clear that if $d\in D_{n}$ then all distances
$k\cdot d$ ($k$ a positive integer) belong to $D_{n}$.
\\
\special{em: graph figure4.gif}
\vskip 1.3truecm
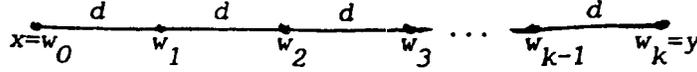
\centerline{Figure 4}
\centerline{$|x-y|=k\cdot d$}
\centerline{$S_{xy}=\bigcup \{S_{ab}:a,b\in \{w_{0},w_{1},...,w_{k}\}, 
|a-b|=d \vee |a-b|=2\cdot d\}$}
\vskip 0.5truecm
 From Figure 5 it is clear that if $d\in D_{n}$, then all distances
$d/k$ ($k$ a positive integer) belong to $D_{n}$. 
Hence $D_{n}\supseteq Q^{+}$; here and subsequently $Q^{+}$ denotes
the set of positive rational numbers.
\\
\special{em: graph figure5.gif}
\vskip 3.3truecm
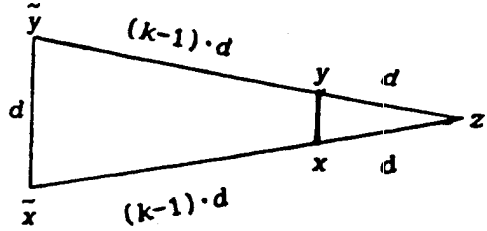
\centerline{Figure 5}
\centerline{$|x-y|=d/k$}

\centerline{$S_{xy}=
S_{\wt{x}\wt{y}}
\cup
S_{\wt{x}x}\cup
S_{xz}\cup
S_{\wt{x}z}
\cup
S_{\wt{y}y}
\cup S_{yz}
\cup
S_{\wt{y}z}$}
\vskip 0.5truecm
{\bf Lemma}. If $x,y\in {\R}^{n}$ ($n>1$) and $\varepsilon >0$
then there exists a finite set $T_{xy}(\varepsilon) \subseteq
{\R}^{n}$ containing $x$ and $y$ such that each map
$f:T_{xy}(\varepsilon )\rightarrow {\R}^{n}$
preserving all unit distances preserves the distance between
$x$ and $y$ to within $\varepsilon $ i.e.
$||f(x)-f(y)|-|x-y||\leq \varepsilon $.
\\
{\it Proof}. It follows from Figure 6.
\\
\special{em: graph figure6.gif}
\vskip 1.8truecm
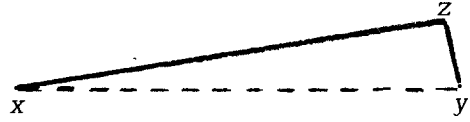
\centerline{Figure 6}
\centerline{$|x-z|,|z-y|\in Q^{+}$, $|z-y|\leq \varepsilon /2$}
\centerline{$T_{xy}(\varepsilon )=S_{xz}\cup S_{zy}$}
\vskip 0.5truecm
{\bf 2. Rigid graphs}

Let $({\R}^{n},1)$ denote the graph with vertex set ${\R}^{n}$ and
with two vertices $u,v$ adjacent if and only if $|u-v|=1$.
A {\it unit-distance graph} in ${\R}^{n}$ is any subgraph of
$({\R}^{n},1)$. A graph $G$ in ${\R}^{n}$ with vertex set $X$
is called {\it rigid} if there exists a $\delta >0$ such that any
$g:X\rightarrow {\R}^{n}$ satisfying:
$|g(x)-x|<\delta $ for all $x\in X$ and $|g(x)-g(y)|=|x-y|$ for all
edges $xy$ of $G$, is an isometry on $X$.
\\
H. Maehara proved that for any real algebraic number $r>0$ there
exists a finite rigid unit-distance graph $G$ in ${\R}^{n}$ ($n>1$)
having two vertices $\alpha $ and $\beta $ at distance $r$, see [7]
for $n=2$ and [7]-[8] for $n>2$. Let us denote by $X$ the set of
vertices of $G$. Let us fix: $\delta>0$ from the definition of
rigidity and $p_{1},...,p_{n},p_{n+1}\in {\R}^{n}$ such that
$|p_{i}-p_{j}|=1$ for $1\leq i<j\leq n+1$. Since
$p_{1},...,p_{n},p_{n+1}$ are affinely independent,
the following mapping
\vskip 0.5 truecm
\centerline{$
{\R}^{n}\ni y
\stackrel{\Phi}{\rightarrow} (|p_{1}-y|,...,|p_{n}-y|,|p_{n+1}-y|)\in
{\R}^{n+1}$}
\vskip 0.5truecm
is injective, see [3] p.127.
\vskip 0.5truecm
$\Phi$ is continuous, $X$ is finite. We can choose $\varepsilon >0$
such that for all $x\in X$, $\wt{x}\in {\R}^{n}$
\vskip 0.5 truecm
\centerline{$max_{1\le i \le n+1} ||p_{i}- \wt{x}|-|p_{i}-x||
\leq \varepsilon $
implies $|\wt{x}-x|<\delta $.}
\vskip 0.5truecm
We prove that $S_{\alpha \beta }:=\bigcup_{x \in X}
\bigcup_{i=1}^{n+1} T_{p_{i}x}(\varepsilon )$
satisfies the condition of the theorem. Let us assume that
$f:S_{\alpha \beta }\rightarrow {\R}^{n}$ preserves the distance $1$.
Hence $f$ preserves distances between $p_{i}$ and $p_{j}$ for
$1\leq i<j\leq n+1$. There is an isometry $I:{\R}^{n}\rightarrow
{\R}^{n}$ such that $f(p_{i})=I(p_{i})$ for $1\leq i \leq n+1$.
Since $S_{\alpha \beta }\supseteq T_{p_{i}x}(\varepsilon )$
($1\leq i\leq n+1$), $f$ preserves distances between $p_{i}$ and
$x\in X$ to within $\varepsilon $. Also $I^{-1}\circ f$ preserves
distances between $p_{i}$ and $x\in X$ to within $\varepsilon $.
Therefore, for all $x\in X$
\vskip 0.5truecm
\centerline{$max_{1 \le i \le n+1} ||(I^{-1}\circ f)
(p_{i})-(I^{-1}\circ f)(x)|-|p_{i}-x||\leq \varepsilon $}
\vskip 0.5truecm
i.e.
\vskip 0.5truecm
\centerline{$max_{1 \le i \le n+1} ||p_{i}-(I^{-1}\circ f)(x)|
-|p_{i}-x|| \leq \varepsilon $}
\vskip 0,5truecm
Hence, for all $x\in X$ $|(I^{-1}\circ f)(x)-x|<\delta $.
If $xy$ is an edge of $G$, then $|x-y|=1$.
Hence $|(I^{-1}\circ f)(x)-(I^{-1}\circ f)(y)|=1=|x-y|$.
By the rigidity of $G$, $I^{-1}\circ f$ is an isometry on $X$.
In particular, $I^{-1}\circ f$ preserves the distance between
$\alpha $ and $\beta $. Hence $f$ preserves the distance between
$\alpha $ and $\beta $, which ends the proof.
\vskip 0.5truecm
{\bf Remark 1.} The lemma from the proof of Theorem 1 implies
the classical Beckman-Quarles theorem.
\vskip 0.5truecm
Let $L=(+,\cdot ,\leq ,0,1)$ denote the language of ordered fields.
We say that $X\subseteq {\R}^{n}$ is definable using $L$ if there is
an $L$-formula $\phi(x_{1},...,x_{n})$ such that
$X=\{(a_{1},...,a_{n})\in {\R}^{n}: \phi (a_{1},...,a_{n})$ $holds\}$
\vskip 0.5truecm
{\bf Remark 2}. Only algebraic distances $|x-y|$ have the property
 from Theorem 1. Indeed, if $x_{1},x_{2}\in {\R}^{n}$ $(n>1)$ and 
$S_{x_{1}x_{2}}=\{x_{1},x_{2},...,x_{m}\}$ then
$\{|x_{1}-x_{2}|\} \subseteq \R$
is definable using $L$ by the following $L$-formula $\psi (z)$:
\\
\centerline{$(0 \leq z) \wedge
\forall x_{11}...\forall x_{1n}... \forall x_{m1}...\forall x_{mn}
( \underbrace{...\wedge (x_{i1}-x_{j1})^2+...+(x_{in}-x_{jn})^2=1
\wedge...}_{(i,j) \in \{(i,j): 1\leq i<j \leq m, |x_{i}-x_{j}|=1\}}$}
$\Rightarrow (x_{11}-x_{21})^2+...+(x_{1n}-x_{2n})^2=z^2)$

 From this, we deduce that $|x_{1}-x_{2}|$ is an
algebraic number defined by some quantifier free $L$-formula equivalent
to $\psi (z)$,
see [4].
\vskip 0.5truecm
{\bf Remark 3}. M. Koshelev proved that Theorem 1 remains valid
for rational distances in two-dimensional Minkowski space-time,
see [6].
\vskip 0.5truecm
{\bf Theorem 2}. If $x,y\in {\F}^{n}$ ($n>1$) then there exists
a finite set $S_{xy}\subseteq {\F}^{n}$ containing $x$ and $y$
such that each map from $S_{xy}$ to ${\R}^{n}$ preserving all
unit distances preserves the distance between $x$ and $y$.
\vskip 0.3truecm
{\it Proof.} Let us denote by $C_{n}$ ($n>1$) the set of all
non-negative constructible numbers $d$ with the following property:

If $x,y\in {\F}^{n}$ and $|x-y|=d$ then there exists a finite set
$S_{xy}\subseteq {\F}^{n}$ such that $x,y\in S_{xy}$ and any map 
$f:S_{xy}\rightarrow {\R}^{n}$ that preserves unit distance also
preserves the distance between $x$ and $y$.

Obviously $0,1\in C_{n}$. Similarly as in the proof of Theorem 1
we can prove that:
\vskip 0.2 truecm
1) $\sqrt{2+2/n}\in C_{n}$,
\vskip 0.2 truecm
2) if $d\in C_{n}$ then all distances $r\cdot d$ ($r$ is positive and
rational) belong to $C_{n}$, therefore all non-negative rational
numbers belong to $C_{n}$,
\vskip 0.2 truecm
3) if $x,y\in {\F}^{n}$ ($n>1$) and $\varepsilon >0$ then there
exists a finite set $T_{xy}(\varepsilon )\subseteq {\F}^{n}$
containing $x$ and $y$ such that each map
$f:T_{xy}(\varepsilon )\rightarrow {\R}^{n}$
preserving all unit distances preserves the distance between
$x$ and $y$ to within $\varepsilon $ i.e.
$||f(x)-f(y)|-|x-y||\leq \varepsilon $.
\vskip 0.4 truecm
If $a,b\in C_{n}$, $a>b>0$ then
$\ \sqrt{a^{2}-b^{2}} \in C_{n}$ (see Figure 7).
\\
\special{em: graph figure7.gif}
\vskip 2.4truecm
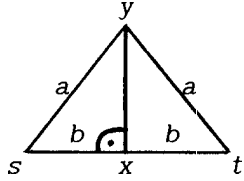
\centerline{Figure 7}
\centerline{$|x-y|$=$\sqrt{a^{2}-b^{2}}$}
\centerline{
$S_{xy}=S_{sx}\cup S_{xt}\cup S_{st}\cup S_{sy}\cup S_{ty}$}
\vskip 0.5truecm
Hence
$\sqrt{3}\cdot a=\sqrt{(2\cdot a)^{2}-a^{2}}\in C_{n}$ and $\sqrt{2}
\cdot a=\sqrt{(\sqrt{3}\cdot a)^{2}-a^{2}}\in C_{n}$.
Therefore $\sqrt{a^{2}+b^{2}}=
\sqrt{(\sqrt{2}\cdot a)^{2}-(\sqrt{a^{2}-b^{2}})^{2}}\in C_{n}$.
\vskip 0.9truecm
Let us see at Figure 8 below, z denotes the centre of the
$(n-1)$-dimensional regular simplex $p_{1}p_{2}...p_{n}$, $n=2$.
This construction shows that if $a,b\in C_{n}$, $a>b>0$, $n\geq 2$
(consistently, equations below Figure 8 describe the general case
$n\geq 2$) then $a-b\in C_{n}$, hence $a+b=2\cdot a-(a-b)\in C_{n}$.
\\
\special{em: graph figure8.gif}
\vskip 4.0truecm
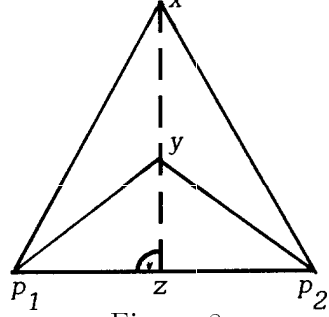
\centerline{Figure 8}
\centerline{$|x-y|=a-b,\hspace{0.2cm} |x-z|=a \in C_{n},\hspace{0.2cm}
|y-z|=b \in C_{n}$}
\centerline{$|p_{i}-p_{j}|=\sqrt{2+2/n}\in C_{n},\hspace{0.2cm}
|z-p_{i}|=\sqrt{1^{2}-(1/n)^{2}}\in C_{n},\hspace{0.2cm}
1 \leq i<j \leq n$}
\centerline{$|x-p_{1}|=\sqrt{|x-z|^{2}+|z-p_{1}|^{2}}=...
=|x-p_{n}|=\sqrt{|x-z|^{2}+|z-p_{n}|^{2}}\in C_{n}$}
\centerline{$|y-p_{1}|=\sqrt{|y-z|^{2}+|z-p_{1}|^{2}}=...
=|y-p_{n}|=\sqrt{|y-z|^{2}+|z-p_{n}|^{2}}\in C_{n}$}
\centerline{
$S_{xy}=\bigcup_{1\le i<j \le n}S_{p_{i}p_{j}}
\cup
\bigcup_{i=1}^{n}S_{xp_{i}}
\cup
\bigcup_{i=1}^{n}S_{yp_{i}}
\cup
T_{xy}(b)$}
\vskip 0.4truecm
In order to prove that $C_{n}\backslash \{0\}$ is a multiplicative
group it remains to observe that if positive
$a,b,c\in C_{n}$ then $\frac{a\cdot b}{c} \in C_{n}$ (see Figure 9).
\\
\special{em: graph figure9.gif}
\vskip 4.5truecm
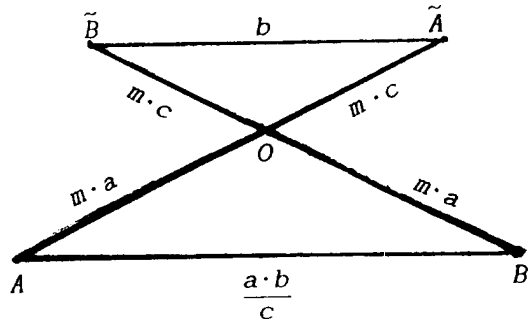
\centerline{Figure 9}
\centerline{$m$ is positive and integer}
\centerline{$b<2\cdot m\cdot c$}
\centerline{$S_{AB}=S_{OA}
\cup
S_{OB}\cup 
S_{O\wt{A}}
\cup
S_{O\wt{B}}
\cup
S_{A\wt{A}}
\cup
S_{B\wt{B}}
\cup
S_{\wt{A}\wt{B}}$}
\vskip 0.6 truecm
If $a\in C_{n}$, $a>1$, then $\sqrt{a}=\frac{1}{2}\cdot 
\sqrt{(a+1)^{2}-(a-1)^{2}} \in C_{n}$;
if $a\in C_{n}$, $0<a<1$,
then $\sqrt{a}=1/\sqrt{\frac{1}{a}}\in C_{n}$.
Thus $C_{n}$ contains all non-negative real numbers contained
in the real quadratic closure of $Q$. This completes the proof.
\vskip 0.5truecm
Directly from the Lemma (see the proof of Theorem 1) we obtain the
following two corollaries, these corollaries and the next elementary
proposition are useful in our next proofs.
\vskip 0.5truecm
{\bf Corollary 1}. If $p\neq q\in {\R}^{n}$ ($n>1$) then any map
$f:T_{pq}(|pq|/2) \rightarrow {\R}^{n}$ that preserves unit distance
satisfies $f(p)\neq f(q)$.
\vskip 0.5truecm
{\bf Corollary 2}. If $J,K,L,M\in {\R}^{n}$ ($n>1$) and $|JK|<|LM|$
then any map
$f:T_{JK}(\frac{|LM|-|JK|}{3})\cup T_{LM}(\frac{|LM|-|JK|}{3})
\rightarrow {\R}^{n}$ that preserves unit distance satisfies
$|f(J)f(K)|<|f(L)f(M)|$.
\vskip 0.5truecm
{\bf Proposition.} If $J,K,L,M\in {\R}^{n}$ ($n>1$) and $|JK|=|LM|$
then there exist $A,B\in {\R}^{n}$ such that:
\\
(1) segments $JK$ and $AB$ are symmetric with respect to some
hyperplane $H(JK,AB)$
\\
and
\\
(2) segments $AB$ and $LM$ are symmetric with respect to some
hyperplane $H(AB,LM)$.
\\
{\it Proof}. Throughout the proof for any two points
$x\in {\R}^{n}$ and $y\in {\R}^{n}$  $H_{xy}$ stands for the
hyperplane which reflects $x$ into $y$, $H_{xy}$ is unique if
$x\neq y$. Let \ $\wt{H}_{xy}:{\R}^{n}\rightarrow {\R}^{n}$ denote
the hyperplane symmetry determined by $H_{xy}$. We put:
$A:=L=\wt{H}_{JL}(J)$, $B:=\wt{H}_{JL}(K)$, $H(JK,AB):=H_{JL}$.
Since $|AB|=|JK|=|LM|=|AM|$ there exists a hyperplane $H_{BM}$
satisfying $\wt{H}_{BM}(A)=A=L$. Therefore, we put $H(AB,LM):=H_{BM}$,
obviously $\wt{H}(AB,LM)(B)=M$. This completes the proof.
\vskip 0.5truecm
Theorems 3 and 4 below yield information about unit-distance
preserving mappings on finite subsets of ${\R}^{n}$.
We present their proofs omitting details.
\vskip 1.5truecm
{\bf Theorem 3.} If $X_{1},X_{2},...,X_{m}\in {\R}^{n}$ ($n>1$)
lie on some affine hyperplane then there exists a finite set
$L(X_{1},X_{2},...,X_{m})\subseteq {\R}^{n}$ containing
$X_{1},X_{2},...,X_{m}$ such that each map from
$L(X_{1},X_{2},...,X_{m})$ to ${\R}^{n}$ preserving all
unit distances preserves the property that
$X_{1},X_{2},...,X_{m}$ lie on some affine hyperplane.
If $X_{1},X_{2},...,X_{m}\in {\F}^{n}$ then we can choose
$L(X_{1},X_{2},...,X_{m})\subseteq {\F}^{n}$.

{\it Proof.} We choose an affine hyperplane $G$ such that
$X_{1},X_{2},...,X_{m} \in G$.
We fix $P\in G$ and choose an integer $u$ such that
$|PX_{i}|\leq 2\cdot u$ for $1\leq i\leq m$. We fix
$O\in {\R}^{n}$ such that $PO$ is perpendicular to $G$
and $|PO|=4\cdot u$. Let us see at Figure 10 below,
this drawing shows the case $n=2$, but equations below
Figure 10 describe the general case $n\geq 2$.
\\
\special{em: graph figure10.gif}
\vskip 9.0truecm
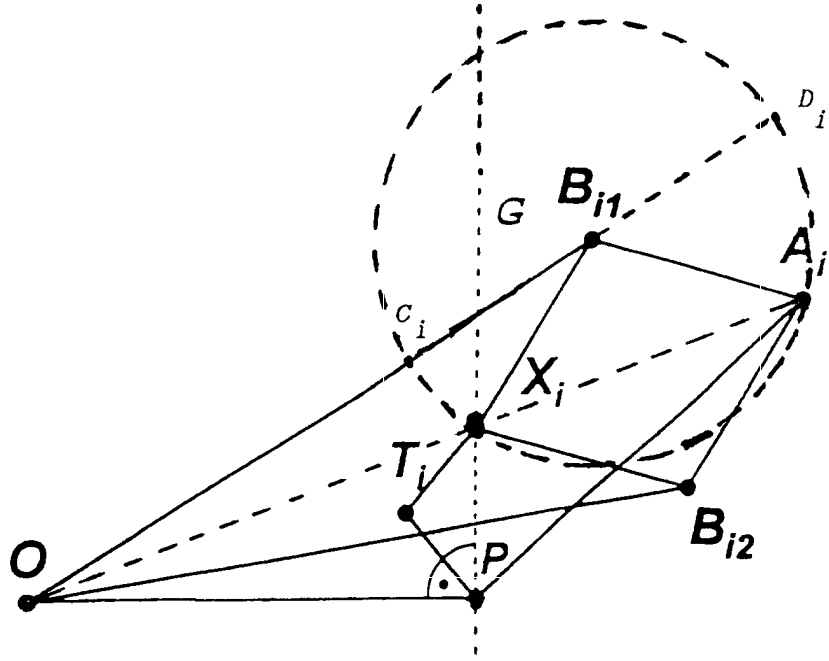
\centerline{Figure 10}
\centerline{$\forall i\in \{1,2,...,m\}\forall j\in \{1,2,...,n\}
|X_{i}B_{ij}|=|A_{i}B_{ij}|=2\cdot u \wedge |OB_{ij}|=6\cdot u$}
\centerline{$|PO|=4\cdot u$, $\forall i\in \{1,2,...,m\}
|PA_{i}|=4\cdot u \wedge |PT_{i}|=|T_{i}X_{i}|=u$}
\vskip 0.35truecm
Erasing all dashed lines and segments $PT_{i}$, $T_{i}X_{i}$, $PA_{i}$,
$PO$ from Figure 10 we obtain the Peaucellier inversor:$|OX_{i}|\cdot
|OA_{i}|=|OC_{i}| \cdot |OD_{i}|=(6\cdot u-2\cdot u)\cdot
(6\cdot u+2\cdot u)=(6\cdot u)^{2}-(2\cdot u)^{2}$
(see [9] in the case $n=2$). Always: $|PX_{i}|\leq 2\cdot u<|PA_{i}|$.
Therefore, if all edge-lengths (continuous segments) are preserved and
if all relations $B_{ik}\neq B_{il}$ ($1\leq k<l\leq n$) are preserved
then $PX_{i}$ is perpendicular to $PO$. From this, we conclude that
\vskip 0.1truecm
\centerline{$L_{i}:=S_{PO}\cup S_{PT_{i}}\cup S_{T_{i}X_{i}}\cup
S_{PA_{i}}\cup ( \bigcup_{j=1}^{n} S_{X_{i}B_{ij}}\cup
S_{A_{i}B_{ij}}\cup S_{OB_{ij}})
\cup$}
\centerline{
$\bigcup \{T_{B_{ik}B_{il}}(|B_{ik}B_{il}|/2):1\leq k<l\leq n\}$}
\vskip 0.1truecm
is adequate for preserving the property that $PX_{i}$ is perpendicular
to $PO$. Hence $L(X_{1},X_{2},...,X_{m}):=\bigcup_{i=1}^{m} L_{i}$
satisfies the condition of the theorem.
This completes the proof for ${\R}^{n}$.
The proof for ${\F}^{n}$ is similar.
\vskip 0.1truecm
{\bf Theorem 4} (cf. Corollary 2). If $J,K,L,M\in {\R}^{n}$ ($n>1$)
and $|JK|=|LM|$ then there exists a finite set $C_{JKLM}\subseteq
{\R}^{n}$ containing $J,K,L,M$ such that any map $f:C_{JKLM}
\rightarrow {\R}^{n}$ that preserves unit distance satisfies
$|f(J)f(K)|=|f(L)f(M)|$. If $J,K,L,M \in {\F}^{n}$ then we can choose
$C_{JKLM} \subseteq {\F}^{n}$.

{\it Proof.} Let us observe that if $J,K,A,B,L,M \in {\R}^{n}$ ($n>1$)
and $|JK|=|AB|=|LM|$ then $C_{JKLM}:=C_{JKAB}\cup C_{ABLM}$ satisfies
the condition of Theorem 4.
By Proposition it is sufficient to prove Theorem 4
in the case where two segments of equal lengths are
symmetric with respect to some hyperplane $H$.

Let us see at figures below, drawings show the case $n=2$, but
equations below figures describe the general case $n\geq 2$.
\\
\special{em: graph figur11a.gif}
\vskip 5.5truecm
\centerline{Figure 11a}
\centerline{$s,t$ are positive and integer}
\centerline{$\forall j\in \{1,2,...,n\}|X_{j}J|=|X_{j}L|=s \wedge
|Y_{j}K|=|Y_{j}M|=t$}
\centerline{$|JK|=|LM|<|JM|=|LK|$ \hspace{0.5cm}
$\delta:=\frac{|JM|-|JK|}{3}$}
\newpage
\special{em:graph figur11b.gif}
\centerline{}
\vskip 5.5truecm
\centerline{Figure 11b}
\centerline{$s,t$ are positive and integer}
\centerline{$\forall j\in \{1,2,...,n\}
|X_{j}J|=|X_{j}L|=s \wedge |Y_{j}K|=|Y_{j}M|=t$}
\centerline{$|JK|=|LM|>|JM|=|LK|$
\hspace{0.5cm}
$\delta :=\frac{|JK|-|JM|}{3}$}
If all edge-lengths (continuous segments) are preserved and if all
relations: $J\neq L$,$X_{k}\neq X_{l}$,$K\neq M$,$Y_{k}\neq Y_{l}$
($1\leq k<l\leq n)$ are preserved then images of $J$ and $L$ are
symmetric and images of $K$ and $M$ are symmetric. Hence:
\\
\begin{math}
C_{JKLM}:=L(X_{1},X_{2},...,X_{n},Y_{1},Y_{2},...,Y_{n}) \cup
(\bigcup_{j=1}^{n} S_{X_{j}J}\cup S_{X_{j}L})
\cup T_{JL}(|JL|/2) \cup
\end{math}
\begin{math}
\bigcup \{ T_{X_{k}X_{l}}(|X_{k}X_{l}|/2):1\leq k<l\leq n\}
\cup
(\bigcup_{j=1}^{n} S_{Y_{j}K}\cup S_{Y_{j}M}) \cup T_{KM}(|KM|/2)
\cup 
\bigcup \{T_{Y_{k}Y_{l}}(|Y_{k}Y_{l}|/2):1\leq k<l\leq n\}
\cup T_{JK}(\delta)
\cup T_{LM}(\delta)
\cup T_{JM}(\delta)
\cup T_{LK}(\delta)
\end{math}
\\
satisfies the condition of the theorem.

Separate formulas for $C_{JKLM}$ are necessary if $J=L$ or $K=M$.
If $J=L$ and $K=M$ then $C_{JKLM}:=\{J,K\}$, if $J\neq L$ and $K=M$
then we see at Figure 11c, Figure 11c shows the case $n=2$,
but equations below Figure 11c describe the general case $n\geq 2$.
\\
\special{em: graph figur11c.gif}
\vskip 2.5truecm
\centerline{Figure 11c}
\centerline{
$s$ is positive and integer, \hspace{0.5cm}
$\forall j\in \{1,2,...,n\} |X_{j}J|=|X_{j}L|=s$}

$C_{JKLM}=L(X_{1},X_{2},...,X_{n},K)
\cup 
(\bigcup_{j=1}^{n} S_{X_{j}J}
\cup
S_{X_{j}L})
\cup T_{JL}(|JL|/2)
\cup$
\\
\centerline{
$\bigcup \{T_{X_{k}X_{l}}(|X_{k}X_{l}|/2):1\leq k<l\leq n\}$}
\\
Symmetric construction applies to the case: $J=L$ and $K\neq M$.
All cases are explained and the proof for ${\R}^{n}$ is complete.
The proof for ${\F}^{n}$ is similar.
\vskip 0.5cm
{\bf Conjecture.}
Let us assume that
$X\subseteq \underbrace{{\R}^{n}\times...\times {\R}^{n}}_{m-times}
\cong {\R}^{mn}$ ($n>1, m \geq 1$) is definable using $L$ and the
relation $X$ is preserved under isometries of ${\R}^{n}$.
We suppose that if $(x_{1},...,x_{m})\in X$ then there exists
a finite set $ P(x_{1},...,x_{m}) \subseteq {\R}^{n}$ containing
$x_{1},...,x_{m}$ such that any map
$f:P(x_{1},...,x_{m}) \rightarrow {\R}^{n}$ that preserves unit
distance satisfies $(f(x_{1}),...,f(x_{m}))\in X$.
\\
\vskip 0.5 truecm
\centerline{{\bf References}}
\par
\begin{enumerate}
\item F. S. Beckman and D. A. Quarles Jr.,
On isometries of euclidean spaces,
{\it Proc. Amer. Math. Soc.}, 4 (1953), 810-815.
\item W. Benz, {\it Geometrische Transformationen (unter besonderer
Ber\"uck\-sichtigung der Lorentztransformationen)}.
BI Wissenschaftsverlag, Ma\-nnheim, Leipzig, Wien,
Z\"urich, 1992.
\item K. Borsuk, {\it Multidimensional analytic geometry}, Polish
Scientific Publishers, Warsaw, 1969.
\item C. C. Chang and H. J. Keisler, {\it Model Theory}, 3rd ed.,
North-Holland, Amsterdam, 1990.
\item U. Everling, Solution of the isometry problem stated
by K. Ciesielski,
{\it Math. Intelligencer}, 10 (1988), No.4, p.47.
\item M. Koshelev, We can measure any distance or any amount
of time with a most primitive clock and a most primitive ruler:
a space-time version of Tyszka's result, {\it Geombinatorics},
7 (1998), No.3, 95-100.
\item H. Maehara, Distances in a rigid unit-distance graph in
the plane, {\it Discrete Appl. Math.} 31 (1991), 193-200.
\item H. Maehara, Extending a flexible unit-bar framework to a
rigid one, {\it Discrete Math.} 108 (1992), 167-174.
\item H. Rademacher and O. Toeplitz,
{\it The enjoyment of mathematics},
Princeton University Press, Princeton, 1994.
\item A. Tyszka, A discrete form of the Beckman-Quarles theorem,
{\it Amer. Math. Monthly} 104 (1997), 757-761.
\end{enumerate}
\vskip 0.2truecm
{\it Technical Faculty}
\\
{\it Hugo Ko{\l}{\l}\c{a}taj University}
\\
{\it Balicka 104, PL-30-149 Krak\'ow, Poland}
\\
{\it rttyszka@cyf-kr.edu.pl}
\\
{\it http://www.cyf-kr.edu.pl/\symbol{126}rttyszka}
\end{document}